\documentclass[11pt]{article}
\usepackage[T1]{fontenc}
\usepackage{lmodern,amsmath,amssymb,amsthm,mathtools,booktabs,array,longtable}
\usepackage[margin=1in]{geometry}
\usepackage{microtype,tikz,float}
\usepackage[colorlinks=true,linkcolor=blue,citecolor=blue,urlcolor=blue]{hyperref}
\hypersetup{pdftitle={A New Record for the Katz–Tao Sum–Difference Exponent},pdfauthor={Atrin Farnamfar},pdfsubject={Additive combinatorics; entropy inequalities; computer-assisted proof}}
\newtheorem{theorem}{Theorem}
\newtheorem{proposition}[theorem]{Proposition}
\newtheorem{lemma}[theorem]{Lemma}
\newtheorem{corollary}[theorem]{Corollary}
\theoremstyle{definition}\newtheorem{definition}[theorem]{Definition}
\theoremstyle{remark}
\newcommand{\cl}{\operatorname{cl}}
\newcommand{\HH}{\mathsf H}
\newcommand{\Q}{\mathbb Q}\newcommand{\Z}{\mathbb Z}
\newcommand{\A}{\mathcal A}\newcommand{\Rel}{\mathcal R}
\newcommand{\Cone}{C_{\!*}}
\title{A New Record for the Katz--Tao\\Sum--Difference Exponent}
\author{Atrin Farnamfar}
\date{October 27, 2026}
\begin{document}
\maketitle
\begin{abstract}
We prove the record bound
\[
C_{3b}\le\frac{209357544569848128}{114631769327023765}
=1.8263483657186\ldots
\]
for the Katz--Tao sum--difference exponent, valid in every abelian group without $2$-torsion. Using only Shannon inequalities, conditional copies and Katz--Tao's original directions, whose classical single-grid argument gives $11/6$, we break $11/6$, obtaining $2762/1509$ and $152/83$. For every abelian group we prove the universal bound
\[
C_{\mathrm{all}}=\frac{10304335425987192362}{5641833174538280329}
=1.8264161855212\ldots,
\]
which is below the previous record for the real Katz--Tao constant. In the AlphaEvolve study~\cite{GGTW}, the search for extremal distributions improved this constant's lower bound in the eighth decimal place and left its upper bound at $11/6$. Our upper bounds are proofs discovered by CASPAC, the computational mathematics engine of the CASP research system, and certified by independent exact checkers. The proofs use stacked entropy certificates; we supply the complete rational identities, proofs of every auxiliary inequality, and exact verification over the stated group classes.
\end{abstract}

\section{Introduction and main results}
The Katz--Tao sum--difference problem asks how many distinct differences a family of pairs can have when its two coordinate projections and its sum projection are small. Sum--difference exponents are a basic ingredient of the arithmetic approach to the Kakeya problem~\cite{KT,GGTW}. We lower the best known upper bound for this sum--difference exponent to $1.8263483657186\ldots$. We also prove the bound $1.8264161855212\ldots$ uniformly in every abelian group, placing the universal bound below the previous record.

For a finite nonempty set $E\subseteq G^2$, where $G$ is an abelian group, write
\[
\pi_x(E)=\{x:(x,y)\in E\},\quad \pi_y(E)=\{y:(x,y)\in E\},\quad
\pi_\pm(E)=\{x\pm y:(x,y)\in E\}.
\]
The problem is to find a small exponent $C$ such that
\begin{equation}\label{eq:projection}
|\pi_-(E)|\le M^C,\qquad
M=\max\{|\pi_x(E)|,|\pi_y(E)|,|\pi_+(E)|\}.
\end{equation}
Over the reals the least exponent is denoted $C_{3b}$. Katz and Tao established $11/6$ in 1999 \cite{KT}. The lower bound and successive upper-bound records are
\begin{center}
\begin{tabular}{@{}llll@{}}\toprule
Side & Work & Year & Bound for $C_{3b}$\\\midrule
Lower (rounded) & Lemm; AlphaEvolve~\cite{Lemm,GGTW} & 2015/2025 & $\approx1.779$\\\midrule
Upper & Katz--Tao & 1999 & $11/6=1.8333333333\ldots$\\
Upper & Tomar~\cite{Tomar} & 2026 & $1.8272962152674\ldots$\\
\textbf{Upper} & \textbf{This work} & \textbf{2026} & $\boldsymbol{1.8263483657186\ldots}$\\\bottomrule
\end{tabular}
\end{center}

For every abelian group the previous best bound was $9213/5036=1.8294281\ldots$~\cite{Tomar}; Theorem~\ref{thm:all} lowers it to $1.8264161855\ldots$, a decrease exceeding $0.0030$.

The constant $C_{3b}=C(\{0,1,\infty\};-1)$ is featured in Problem~6.30 (arithmetic Kakeya conjecture) of the AlphaEvolve study by Georgiev, G\'omez-Serrano, Tao and Wagner~\cite{GGTW}. There, AI search produced extremal distributions that improved the lower bound in the eighth decimal place, while the reported upper bound remained $11/6$. Such distributions witness lower bounds for this constant; a universal upper bound requires a proof. Here, CASPAC produces exactly verifiable proof certificates, establishing new records on the upper-bound side of the same problem.

The classical single-grid entropy argument gives $11/6$ using Shannon inequalities, conditional copies and the original directions; see~\cite{KT} and~\cite[Section~3 and Remark~3.4]{Tomar}. With those same ingredients, we obtain $2762/1509=11/6-3/1006$ and $152/83=11/6-1/498$. Thus $11/6$ is not the limit of the original-direction method: the classical proof can be improved using its own tools. These bounds hold in every abelian group and use neither scalar cross forms nor auxiliary non-Shannon inequalities.

\paragraph{How the certificates were found.}
CASP is a research system that carries a long investigation across many AI sessions, retaining what was tried, what failed and what worked. CASPAC (CASP Arith Compute) is its computational mathematics engine, developed by the author. Within this engine, AI agents explored stacked configurations, discovered the certificates and verified the resulting proofs with independent exact checkers, under the author's direction. CASP used Anthropic Claude and OpenAI Codex models. The resulting proof objects can be verified independently of the system that found them.

The new ingredient throughout is stacking: complete pages are copied conditionally on chosen readings and glued into one configuration, with equal laws carried across the copies (Lemma~\ref{lem:enriched}). Stacking breaks $11/6$ with the original directions alone (Theorem~\ref{thm:standard}) and yields the records of Theorems~\ref{thm:record} and~\ref{thm:all}.

All random variables in this paper have finite support, and the entropy base is arbitrary and fixed. A group is \emph{without $2$-torsion} when $2g=0$ implies $g=0$; this includes torsion-free groups and finite abelian groups of odd order.

\begin{theorem}[Record bound for the Katz--Tao exponent]\label{thm:record}
Let $X,Y$ take values in an abelian group without $2$-torsion. Then
\begin{equation}\label{eq:weighted}
\HH(X-Y)\le a\HH(X)+b\HH(Y)+c\HH(X+Y),
\end{equation}
where
\begin{align*}
a&=\frac{80908638342740092}{114631769327023765},\\
b&=\frac{74518340847137328}{114631769327023765},\\
c&=\frac{53930565379970708}{114631769327023765}.
\end{align*}
Consequently \eqref{eq:projection} holds with
\begin{equation}\label{eq:record}
\Cone=a+b+c=\frac{209357544569848128}{114631769327023765}
=1.8263483657186\ldots.
\end{equation}
In particular $C_{3b}\le\Cone$.
\end{theorem}

\begin{theorem}[Record bound for every abelian group]\label{thm:all}
For every finitely supported pair in an abelian group, inequality \eqref{eq:weighted} holds with
\begin{align*}
a&=\frac{31837334702369419045}{45134665396306242632},\\
b&=\frac{29267304885161681931}{45134665396306242632},\\
c&=\frac{21330043820366437920}{45134665396306242632}.
\end{align*}
Thus \eqref{eq:projection} holds with
\begin{equation}\label{eq:all}
C_{\mathrm{all}}=a+b+c=\frac{10304335425987192362}{5641833174538280329}
=1.8264161855212\ldots.
\end{equation}
\end{theorem}

\begin{theorem}[Breaking $11/6$ with the original vertex forms]\label{thm:standard}\label{thm:fivepage}
For every finitely supported pair in an abelian group,
\begin{align}
1509\HH(X-Y)&\le1077\HH(X)+962\HH(Y)+723\HH(X+Y),\label{eq:standardbound}\\
83\HH(X-Y)&\le59\HH(X)+52\HH(Y)+41\HH(X+Y).\label{eq:fivepage}
\end{align}
Both have conditional-copy proofs using only Shannon inequalities and the standard vertex forms
\begin{equation}\label{eq:standard}
x_v,\quad y_v,\quad s_v=x_v+y_v,\quad d_v=x_v-y_v.
\end{equation}
The first inequality gives the bound
\[
C_{\mathrm{std}}=\frac{2762}{1509}=\frac{11}{6}-\frac{3}{1006}=1.8303512259774\ldots.
\]
The second inequality gives the bound
\[
C_5=\frac{152}{83}=\frac{11}{6}-\frac{1}{498}=1.8313253012048\ldots.
\]
\end{theorem}

The entropy tuples may involve many vertices, but every scalar variable in Theorem~\ref{thm:standard} is one of \eqref{eq:standard} at a single vertex.

\begin{table}[ht]
\centering\small
\begin{tabular}{@{}llll@{}}\toprule
Result & Group class & Certificate ingredients & Exponent\\\midrule
Theorem~\ref{thm:record} & No $2$-torsion & Cross forms; proved lemmas & $1.8263483657\ldots$\\
Theorem~\ref{thm:all} & Every abelian group & Cross forms; proved lemmas & $1.8264161855\ldots$\\
\eqref{eq:standardbound} & Every abelian group & Standard forms; Shannon & $2762/1509$\\
\eqref{eq:fivepage} & Every abelian group & Standard forms; Shannon & $152/83$\\
\bottomrule
\end{tabular}
\caption{The four exact bounds and their proof ingredients.}
\label{tab:results}
\end{table}

Section~\ref{sec:configuration} specifies the literal certificate configurations. Sections~\ref{sec:closure}--\ref{sec:certificate} establish their group validity and exact coefficient identities.

\section{Entropy and conditional copies}\label{sec:entropy}
All random variables below have finite support. We use Shannon entropy in any fixed base greater than one and write $\HH[U]$ for the joint entropy of the variables in a finite set $U$. Repeated variables do not change this entropy. We use the standard facts
\begin{align}
\HH[U]+\HH[V]&\ge\HH[U\cup V]+\HH[U\cap V],\label{eq:sub}\\
\HH[U]&\le\HH[V]\quad(U\subseteq V),\label{eq:mono}\\
\HH[U\cup V]&=\HH[U]\quad\text{if $V$ is determined by $U$}.
\label{eq:det}
\end{align}
Submodularity is equivalent to nonnegativity of conditional mutual information. Equal laws give equal entropies. These facts apply to arbitrary jointly distributed variables.

\begin{lemma}[Conditional-copy extension]\label{lem:copy}
Let $O$ be a finite random tuple, $U$ a subtuple of $O$, and $G=g(U)$ for a fixed function $g$. There is an extension of the probability space carrying $U'$ such that
\begin{enumerate}
\item $(U',G)$ and $(U,G)$ have the same law;
\item $g(U')=g(U)$ almost surely;
\item $U'$ is conditionally independent of $O$ given $G$.
\end{enumerate}
Consequently
\begin{equation}\label{eq:copy}
\HH[O,U']+\HH[G]-\HH[O]-\HH[U']=0.
\end{equation}
\end{lemma}
\begin{proof}
For an old outcome $o$, sample $u'$ with conditional law
$\mathbb P(U=u'\mid G=g(U(o)))$. Equivalently, give $(o,u')$ probability
\[
\mathbb P(O=o)\,\mathbb P(U=u'\mid G=g(U(o))).
\]
Here $U(o)$ denotes the source subtuple of $o$. Values of $G$ with zero probability need no specification. Summing over $o$ shows that the marginal law of $U'$ equals that of $U$ and that the joint laws with $G$ agree. The conditional distribution is supported on $g(u')=g(U(o))$ and depends on $o$ only through $g(U(o))$. This proves the three properties. Conditional independence gives
$\HH[O,U'\mid G]=\HH[O\mid G]+\HH[U'\mid G]$; because $G$ is determined by both $O$ and $U'$, expansion gives \eqref{eq:copy}.
\end{proof}

\begin{definition}
A pair $(X,Y)$ is \emph{difference-injective} if the map $(x,y)\mapsto x-y$ is injective on its support. Then $\HH(X,Y)=\HH(X-Y)$, and the pair is determined by its difference.
\end{definition}
Every copied vertex retains this property because it retains the law of the original pair. We first prove the inequality for difference-injective pairs. The following reduction removes that restriction.

\begin{lemma}[Removal of the injectivity restriction]\label{lem:selector}
Suppose $a,b,c\ge0$ and
\[
\HH(X,Y)\le a\HH(X)+b\HH(Y)+c\HH(X+Y)
\]
for every difference-injective pair in a fixed abelian group. Then
\[
\HH(X-Y)\le a\HH(X)+b\HH(Y)+c\HH(X+Y)
\]
for every finitely supported pair in that group.
\end{lemma}
\begin{proof}
Put $D=X-Y$. Independently of $D$, choose a random selector $\Phi$ which, for each possible difference $d$, selects a pair with the conditional law of $(X,Y)$ given $D=d$. The choices for different $d$ may be made independently. Set $(\widetilde X,\widetilde Y)=\Phi(D)$. This pair has the original law and difference $D$.

For every fixed selector, its selected pairs have distinct differences. Applying the assumed inequality to the conditional law given $\Phi$ and averaging yields
\[
\HH(\widetilde X,\widetilde Y\mid\Phi)
\le a\HH(\widetilde X\mid\Phi)+b\HH(\widetilde Y\mid\Phi)
   +c\HH(\widetilde X+\widetilde Y\mid\Phi).
\]
Given $\Phi$, the pair $\Phi(D)$ and $D$ determine each other. The left side is therefore $\HH(D\mid\Phi)=\HH(D)$. Conditioning does not increase entropy, and the three coefficients are nonnegative. Replacing each term on the right by its unconditional entropy proves the claim.
\end{proof}

\section{Certificate configurations}\label{sec:configuration}
Fix a difference-injective pair $Z=(X,Y)$. Every vertex below is a pair with the law of $Z$. Letters $A,B,C$ mean the readings $x,y,x+y$, respectively; a dot means no reading.

\begin{lemma}[Shared and transported readings]\label{lem:enriched}
Let $U'$ be a conditional copy of a tuple $U$ over $V=g(U)$. Every reading $f(U)$ determined by $V$ satisfies $f(U')=f(U)$ almost surely. More generally, for any finite list of readings $f_1,\ldots,f_m$,
\[
(f_1(U),\ldots,f_m(U),V)\ \stackrel{\mathrm{law}}{=}\
(f_1(U'),\ldots,f_m(U'),V).
\]
\end{lemma}
\begin{proof}
The conditional law is supported on the fibre $g(U')=V$. If $f(U)=\varphi(V)$ on the support, the same relation holds for $U'$ and gives $f(U')=\varphi(V)=f(U)$. Equality of the joint laws of $(U,V)$ and $(U',V)$ then gives the displayed identity after applying the fixed readings. The conditional-copy extension places no further restriction on $g$.
\end{proof}

\subsection{The initial page for the cross-form bounds}
Use the following three steps, with new vertices numbered in creation order:
\begin{enumerate}
\item Start with $z_0=Z$. Copy $z_0$ over $x_0$ to obtain $z_1$.
\item Copy $(z_0,z_1)$ over $(y_0,y_1)$ to obtain $(z_2,z_3)$.
\item Copy $(z_0,z_1)$ over $(s_0,y_1)$ to obtain $(z_4,z_5)$.
\end{enumerate}
In each step the new tuple is independent of the entire old tuple conditional on the indicated readings. This is the six-pair page of~\cite{Tomar}. Its seven coordinate relations give the parametrization
\begin{equation}\label{eq:rootparam}
\begin{array}{lll}
z_0=(a,b),&z_1=(a,c),&z_2=(d,b),\\
z_3=(d,c),&z_4=(e,a+b-e),&z_5=(e,c).
\end{array}
\end{equation}
Here $a,b,c,d,e$ parametrize the coordinates; their joint law is specified by the copy steps. This coordinate notation is local to the page parametrization.

The Katz--Tao diagram has top row $(z_2,z_0,z_4)$ and bottom row $(z_3,z_1,z_5)$. Accordingly, all six-letter glue words for these configurations are read in the order
\begin{equation}\label{eq:rootorder}
(2,0,4,3,1,5).
\end{equation}

\begin{figure}[ht]
\centering
\begin{tikzpicture}[x=2.3cm,y=1.25cm,every node/.style={font=\small}]
\foreach \c/\u/\v in {0/2/3,1/0/1,2/4/5}{
\coordinate (u\c) at (\c,1);\coordinate (v\c) at (\c,0);
\draw (u\c)--node[right]{$A$}(v\c);
\fill (u\c) circle (2pt);\fill (v\c) circle (2pt);
\node[above=6pt] at (u\c){$z_{\u}$};\node[below=6pt] at (v\c){$z_{\v}$};}
\draw (u0)--node[above]{$B$}(u1)--node[above]{$C$}(u2);
\draw (v0)--node[below]{$B$}(v1)--node[below]{$B$}(v2);
\end{tikzpicture}
\caption{The initial six-pair page in glue-word order. Edges express equality of readings. Each new complete page is sampled from the specified page law, conditional on its glue readings.}
\end{figure}
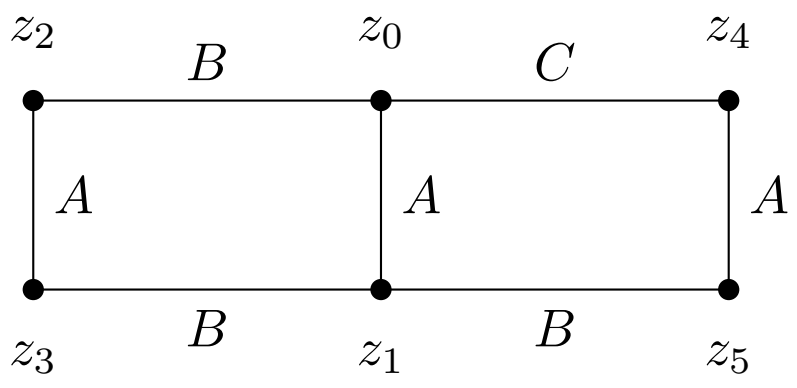

Call this page $P_0$. For a word $w$, let $G_w(P)$ be its tuple of readings in order \eqref{eq:rootorder}. Given all previous pages, sample each new $P_j$ as a conditional copy of the \emph{root page $P_0$} over $G_{w_j}(P_0)$. Thus the pages form a star with source $P_0$ at every step. Table~\ref{tab:glues} gives the complete list of words, in creation order.

\begin{table}[ht]
\centering\small
\begin{tabular}{@{}cll@{}}\toprule
$j$ & Theorem~\ref{thm:record} & Theorem~\ref{thm:all}\\\midrule
1 & \texttt{..BCCA} & \texttt{..B.CC}\\
2 & \texttt{..B.CC} & \texttt{ABAB..}\\
3 & \texttt{ABAB..} & \texttt{..BAAA}\\
4 & \texttt{..BAAA} & \texttt{..CCC.}\\
5 & \texttt{..BACC} & \texttt{..BACC}\\
6 & \texttt{..CCC.} & \texttt{..BCCA}\\\bottomrule
\end{tabular}
\caption{Complete-page star couplings, in creation order, with source $P_0$ at each step.}
\label{tab:glues}
\end{table}

At a page-copy step set $O=(P_0,\ldots,P_{j-1})$, $U=P_0$ and $G=G_{w_j}(P_0)$ in Lemma~\ref{lem:copy}. Every new page has the root-page law and satisfies
\begin{equation}\label{eq:pageci}
\HH[O,P_j]+\HH[G]-\HH[O]-\HH[P_j]=0.
\end{equation}
The old tuple $O$ in every identity is the entire prefix of previously added pages.

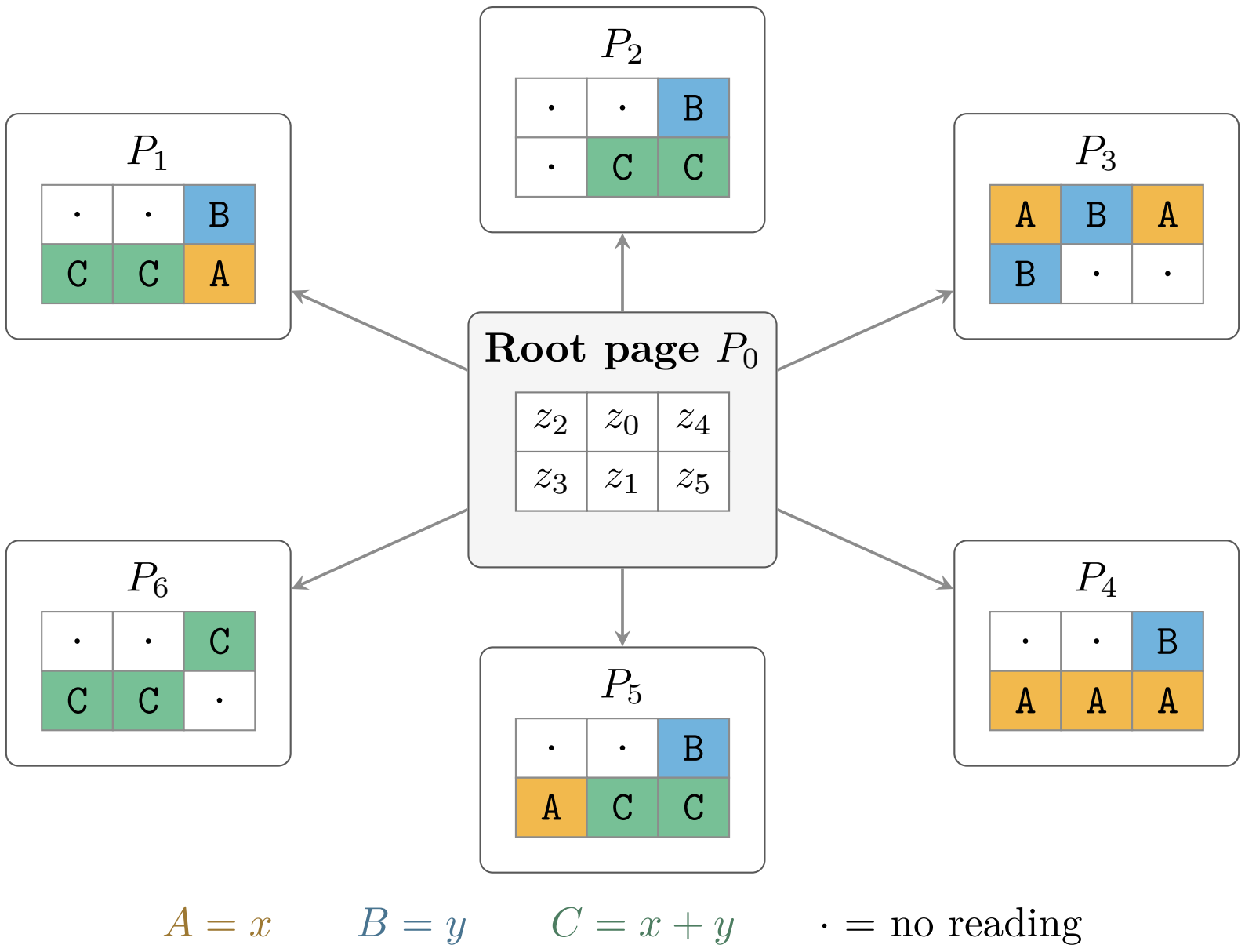
\begin{figure}[htbp]
\centering
\begin{tikzpicture}[x=1cm,y=1cm,>=stealth,every node/.style={font=\small}]
\definecolor{glueA}{RGB}{242,185,77}
\definecolor{glueB}{RGB}{113,179,221}
\definecolor{glueC}{RGB}{119,192,150}
\tikzset{page/.style={draw=black!65,rounded corners=3pt,fill=white,minimum width=2.4cm,minimum height=1.9cm},link/.style={->,draw=black!45,line width=.7pt}}
\node[page,fill=black!4,minimum width=2.6cm,minimum height=2.15cm] (root) at (0,0) {};
\foreach \j/\xx/\yy in {1/-4/1.8,2/0/2.7,3/4/1.8,4/4/-1.8,5/0/-2.7,6/-4/-1.8}{
  \node[page] (p\j) at (\xx,\yy) {};
  \draw[link] (root)--(p\j);
  \node at (\xx,\yy+.62) {$P_{\j}$};
}
\node[font=\bfseries\small] at (0,.75) {Root page $P_0$};
\foreach \u/\v/\label in {0/0/2,1/0/0,2/0/4,0/1/3,1/1/1,2/1/5}{
  \draw[fill=white,draw=black!50] (-.9+.6*\u,.4-.5*\v) rectangle ++(.6,-.5);
  \node at (-.6+.6*\u,.15-.5*\v) {$z_{\label}$};
}
\foreach \xx/\yy/\word in {-4/1.8/{.,.,B,C,C,A},0/2.7/{.,.,B,.,C,C},4/1.8/{A,B,A,B,.,.},4/-1.8/{.,.,B,A,A,A},0/-2.7/{.,.,B,A,C,C},-4/-1.8/{.,.,C,C,C,.}}{
  \foreach \letter [count=\k from 0] in \word {
    \pgfmathtruncatemacro{\col}{mod(\k,3)}
    \pgfmathtruncatemacro{\row}{floor(\k/3)}
    \def\cellcolour{white}
    \def\la{A}\def\lb{B}\def\lc{C}
    \ifx\letter\la\def\cellcolour{glueA}\fi
    \ifx\letter\lb\def\cellcolour{glueB}\fi
    \ifx\letter\lc\def\cellcolour{glueC}\fi
    \draw[fill=\cellcolour,draw=black!45] (\xx-.9+.6*\col,\yy+.35-.5*\row) rectangle ++(.6,-.5);
    \node[font=\small\ttfamily] at (\xx-.6+.6*\col,\yy+.10-.5*\row) {\letter};
  }
}
\node[font=\small] at (0,-4.1) {\textcolor{glueA!65!black}{$A=x$}\qquad\textcolor{glueB!65!black}{$B=y$}\qquad\textcolor{glueC!65!black}{$C=x+y$}\qquad$\cdot$ = no reading};
\end{tikzpicture}
\caption{The seven-page configuration proving Theorem~\ref{thm:record}. The six cells on each outer page display its glue word in order $(2,0,4,3,1,5)$. A coloured cell shares the indicated reading with the corresponding cell of $P_0$. Pages are added in numerical order; the conditional independence at each step is with respect to all earlier pages.}
\label{fig:stack}
\end{figure}

\subsection{Cross forms and auxiliary information inequalities}
Both record certificates contain the 43 root-page cross forms of~\cite{Tomar}. Theorem~\ref{thm:record} adds 35 explicitly listed forms on copied pages, for 78 extra scalar forms in total; Theorem~\ref{thm:all} adds 42, for a total of 85. Appendix~\ref{app:forms} defines the root-form indexing and lists every transported form.

Theorem~\ref{thm:record} uses the 72-lemma library of~\cite{Tomar}. It consists of the Mat\'u\v{s} family with parameters $1,\ldots,62$ and ten Dougherty--Freiling--Zeger inequalities~\cite{Matus,DFZ}. Theorem~\ref{thm:all} extends this library to 77 entries with five further instances of the first Dougherty--Freiling--Zeger extension rule~\cite[Section~IX]{DFZ}, named \texttt{dfz109\_R1\_4} through \texttt{dfz109\_R1\_8}. Each carries its complete conditional-copy proof, checked by the same verifier. The two main certificates invoke twelve and eleven lemma rows, respectively.

Every library entry includes its copy steps and a proof from Shannon inequalities, copy identities and equal laws, possibly applying earlier verified entries. The verifier checks this dependency graph and each coefficient identity before any application in the main certificate. Every auxiliary inequality is therefore verified from its complete conditional-copy proof. A library copy may condition on any old tuple $V$: apply Lemma~\ref{lem:copy} to the joint source $(U,V)$ over $V$, and retain the new $U$ component. This gives the usual general conditional-copy rule used in those proofs.

\subsection{Chains with the original vertex forms}
For \eqref{eq:standardbound}, build a six-pair page with positions $(r,c)$, $r\in\{0,1\}$ and $c\in\{0,1,2\}$, as follows. Copy the first vertex over $x$ to create the first column. Copy that column over its two $y$ readings to create the second column. Copy the \emph{second} column over its top $s$ and bottom $y$ readings to create the third. The top edges have colours $B,C$, the bottom edges $B,B$, and all three vertical edges have colour $A$.

The pages form a chain: for $j=1,\ldots,17$, copy the preceding page $P_{j-1}$ over $G_{w_j}(P_{j-1})$, independently of the entire old prefix conditional on this tuple. Read the words in row-major order
\[
(0,0),(0,1),(0,2),(1,0),(1,1),(1,2).
\]
The seventeen words, in order, are
\begin{equation}\label{eq:chainwords}
\begin{gathered}
\texttt{..CCC.},\ \texttt{..B.CC},\ \texttt{..BCA.},\ \texttt{..BACC},\ \texttt{..BCCA},\\
\texttt{.B...A},\ \texttt{.B.ABA},\ \texttt{AABB..},\ \texttt{..BC.C},\ \texttt{.BBA.B},\\
\texttt{...A.A},\ \texttt{B...C.},\ \texttt{B..CC.},\ \texttt{B.A...},\ \texttt{C..C..},\\
\texttt{C.B..C},\ \texttt{C.BC..}.
\end{gathered}
\end{equation}
This gives 108 vertices and 20 conditional-copy steps, including the three steps that create the first page. Its ground set consists exactly of the four forms \eqref{eq:standard} at each vertex. There are no cross forms and no auxiliary lemma rows.

For \eqref{eq:fivepage}, use the same initial page and chain rule with four glue words:
\begin{equation}\label{eq:fivewords}
\texttt{..CCC.},\qquad\texttt{..BC.C},\qquad
\texttt{C..C..},\qquad\texttt{C.B..C}.
\end{equation}
This gives five pages, 30 vertices and seven conditional-copy steps. It uses the same standard forms and no cross forms or auxiliary lemma rows.

\begin{proposition}[Realization]\label{prop:realization}
Each of the configurations just specified exists for every finitely supported difference-injective pair. Every vertex has that pair's law. The stated linear relations, full-prefix copy identities, and equal-law identifications in the corresponding certificate hold almost surely or as entropy identities, respectively.
\end{proposition}
\begin{proof}
At every step the conditioning tuple is a function of the source tuple. Lemma~\ref{lem:copy} therefore extends the existing law without changing any old marginal, gives the new tuple its source law, and fixes the conditioning readings. Relations internal to the source are inherited by the copy. Induction proves the marginal and relation statements, and \eqref{eq:copy} gives each full-prefix identity. Equality of the joint laws of the source and copy with the conditioning tuple gives equality for every tuple of their corresponding readings. Each vertex retains difference injectivity because it retains the law of $Z$.
\end{proof}

\section{Every abelian group: integrality of the certificates}\label{sec:closure}
The group hypotheses are established by deterministic closure over the integers or the integers with $2$ inverted.
Let $\A$ be the labelled set of standard forms at every vertex, together with the cross forms when present. For a ring $R\in\{\Z,\Z[1/2],\Q\}$, let $\Rel_R$ be the module of relations generated by the copy configuration. Starting from $S\subseteq\A$, repeatedly:
\begin{enumerate}
\item add each labelled form in $\operatorname{span}_R(S)+\Rel_R$;
\item whenever $d_v$ is determined, add $x_v$ and $y_v$.
\end{enumerate}
Denote the resulting finite closure by $\cl_R S$. The second operation uses difference injectivity; it is a deterministic implication, not a linear relation between $x_v,y_v,d_v$.

The relation module is generated step by step. For a source tuple $U$ copied over readings $V$, take the module spanned by the old relations and $\cl_R V$, intersect it with the coordinate forms supported on $U$, and, for every generator $\alpha$ of the intersection, add
\begin{equation}\label{eq:relcopy}
\alpha(U')-\alpha(U)=0.
\end{equation}
The old relation module is retained. This operation inherits internal source relations and identifies the source readings determined by $V$. Intersection is computed by eliminating coordinates outside $U$ first. In the integer calculation all row operations are integral gcd operations; in the localized calculation only powers of $2$ may additionally be inverted.

\begin{lemma}[Soundness over the specified ring]\label{lem:closure}
The generated relations and closures are sound over every abelian group when $R=\Z$, and over every abelian group without $2$-torsion when $R=\Z[1/2]$. In particular,
\begin{equation}\label{eq:closure}
\HH[\cl_R S]=\HH[S].
\end{equation}
\end{lemma}
\begin{proof}
For $R=\Z$, membership of $f$ in the module generated by known forms and valid relations expresses $f$ as an integer combination of known readings, modulo relations that vanish almost surely. Thus $f$ is determined. For $R=\Z[1/2]$, clear denominators to write
\[
2^k f=\sum_{g\in S}n_g g+\rho,
\]
where $n_g$ are integers and $\rho$ is an integer combination of valid relation generators. Multiplication by $2^k$ is injective in the stated group class, so the known readings again determine $f$. Division need not be an operation defined on the whole group.

The difference-injectivity rule is sound at every vertex. At a copy step every $\alpha$ in the intersection defining \eqref{eq:relcopy} is determined by the conditioning tuple and is supported on the source. Equality of the joint laws with that tuple forces the same value in the copy. Applying this argument inductively proves the validity of every generated relation. Iterating the two closure rules proves that the closure is determined by $S$; the reverse determination holds because it contains $S$.
\end{proof}

\subsection{Exact agreement of closures}
The certificates use rational closures as canonical names for finite sets of scalar random variables. For every tuple $T$, define
\[
h(T)=\HH[\cl_{\Q}T].
\]
For every tuple queried in a certificate row, equal-law identification or target term, we verify
\begin{equation}\label{eq:audit}
\cl_R T=\cl_{\Q}T.
\end{equation}
Lemma~\ref{lem:closure} then gives $h(T)=\HH[T]$, so each Shannon inequality, copy identity, equal law and auxiliary information inequality applies directly in the stated group class.

\begin{proposition}[Domain validity]\label{prop:domain}
For the certificate of Theorem~\ref{thm:record}, equality \eqref{eq:audit} holds on every queried tuple with $R=\Z[1/2]$. For Theorem~\ref{thm:all} and both standard-form certificates it holds with $R=\Z$. All conditioning closures, shared-form determinations and relation-space comparisons hold at every copy prefix. Every certificate row is therefore valid over its stated group class.
\end{proposition}
\begin{proof}[Exact finite verification]
The domain verifier generates the relation modules and deterministic closures over the indicated ring. At each copy prefix it checks the conditioning closure and both inclusions between the rational spans of the independently generated relation spaces. For a cross form kept fixed in an equal-law entry, it also checks determination by the glue at that prefix. For transported forms it checks the literal coefficient image under the source-to-copy map. It then checks \eqref{eq:audit} on every tuple queried by the actual proof rows. All arithmetic is exact. The supplement supplies the complete enumerations, passing receipts and implementation; Appendix~\ref{app:reproduction} records their sizes and commands.
\end{proof}

\section{The certificate identities}\label{sec:certificate}
With $h(T)=\HH[\cl_{\Q}T]$ as above, the elemental certificate rows are:
\begin{align*}
\text{sub:}\quad &h(K,i)+h(K,j)-h(K,i,j)-h(K)\ge0,\\
\text{mono:}\quad &h(K,i)-h(K)\ge0,\\
\text{ci:}\quad &h(V,U')+h(O)-h(O,U')-h(V)=0.
\end{align*}
The last row uses the entire old tuple $O$ and is the negative of \eqref{eq:copy}. Its multiplier may have either sign. Equal-law entries identify the entropies of a source tuple and its corresponding copied tuple. The record certificates also use cross-form generators: a reading determined by the glue is kept fixed, and a source-supported reading may be replaced by its literal image on the copy. Lemma~\ref{lem:enriched} proves these equal laws, and Proposition~\ref{prop:domain} verifies the associated closure substitutions.

A lemma row is a verified information inequality with each abstract variable replaced by a tuple of ground-set forms, multiplied by a nonnegative rational. Shannon and monotonicity multipliers are also nonnegative. All lemma proofs are checked before use, including entries not invoked by a main row.

\begin{proposition}[Exact identities]\label{prop:identities}
For each of the four certificates, summing its rows with their listed rational multipliers gives, modulo the listed equal-law identifications,
\begin{equation}\label{eq:identity}
\sum_r m_r R_r
=a\,h(x_0)+b\,h(y_0)+c\,h(s_0)-h(x_0,y_0),
\end{equation}
with the corresponding weights of Theorems~\ref{thm:record}--\ref{thm:standard}. The coefficient residual on every entropy class is exactly zero. Table~\ref{tab:counts} gives the sizes of these finite identities.
\end{proposition}
\begin{proof}[Exact finite verification]
The independent rational checker reconstructs the configuration, generates the closures, and identifies equal-law classes. For every row it reconstructs the entropy terms from the stated row rule and adds the exact rational coefficients. It then subtracts the right side of \eqref{eq:identity}. Every coefficient vanishes. Auxiliary lemma proofs are themselves verified by exact coefficient identities in their abstract conditional-copy configurations. A second rational replay uses the public checker of~\cite{Tomar} with sparse exact closure and explicit support for the shared and transported readings of Lemma~\ref{lem:enriched}. Its Shannon-row and auxiliary-lemma verification routines are unchanged. The independently implemented checker and this extended replay return the same weights and constants. The full verifier code and observed receipts are supplied with the paper.
\end{proof}

\begin{proof}[Proof of Theorems~\ref{thm:record}--\ref{thm:standard}]
First take a difference-injective pair. Proposition~\ref{prop:realization} supplies the joint law; Proposition~\ref{prop:domain} justifies every identity, lemma application and Shannon row in the stated group class. In \eqref{eq:identity}, all inequality rows are nonnegative and all copy identities are zero. Hence
\[
\HH(X,Y)\le a\HH(X)+b\HH(Y)+c\HH(X+Y).
\]
Lemma~\ref{lem:selector} removes the injectivity restriction and gives the asserted entropy inequalities. The configurations and their row lists establish the stated proof ingredients. The projection conclusion follows from Corollary~\ref{cor:projection} below.
\end{proof}

\section{Consequences for projections and progressions}\label{sec:consequences}
\begin{corollary}[Weighted projection inequalities]\label{cor:projection}
For each triple $(a,b,c)$ proved above, in its stated group class,
\begin{equation}\label{eq:product}
|\pi_-(E)|\le |\pi_x(E)|^a|\pi_y(E)|^b|\pi_+(E)|^c
\end{equation}
for every finite nonempty $E\subseteq G^2$. In particular \eqref{eq:projection} holds with $C=a+b+c$.
\end{corollary}
\begin{proof}
Select one pair from $E$ for each distinct difference and put the uniform law on the selected set. This pair is difference-injective and its joint entropy is $\log|\pi_-(E)|$. The entropies of its three input projections are at most the logarithms of the corresponding projection sizes of $E$. Apply the injective entropy inequality and exponentiate.
\end{proof}
For example, \eqref{eq:fivepage} gives the exact integer-power inequality
\[
|\pi_-(E)|^{83}\le
|\pi_x(E)|^{59}|\pi_y(E)|^{52}|\pi_+(E)|^{41}.
\]
The asymmetric weights are useful when the projection sizes differ. Exchanging $X,Y$ gives another valid inequality because negation preserves entropy. Taking the smaller of the two yields
\begin{equation}\label{eq:symmetric}
\HH(X-Y)\le\frac{a+b}{2}\bigl(\HH(X)+\HH(Y)\bigr)
+c\HH(X+Y)-\frac{|a-b|}{2}|\HH(X)-\HH(Y)|.
\end{equation}
This refinement requires no additional certificate.

\begin{corollary}[Distinct progression differences]\label{cor:ap}
Let $A\subset\mathbb R$ be finite with $|A|=n\ge1$. The number of positive $d$ for which $A$ contains a three-term arithmetic progression $a,a+d,a+2d$ is at most $n^{\Cone}$.
\end{corollary}
\begin{proof}
Choose one such progression for each distinct positive $d$ and put the pair $(a,a+2d)$ in $E$. Then $|\pi_x(E)|,|\pi_y(E)|\le n$, and $\pi_+(E)\subseteq\{2a:a\in A\}$, so its size is at most $n$. The differences of these pairs are the distinct values $-2d$. Apply \eqref{eq:projection}; the empty case is immediate.
\end{proof}

\section{Exact verification}
Each bound is a finite rational identity. The specified joint probability law, deterministic relations and valid information inequalities establish the meaning of every row; the checkers establish the exact cancellation of every entropy coefficient. The standalone supplement includes the complete proof objects, all auxiliary inequality proofs, independent verification code and deliberately corrupted inputs used to test rejection. Every check uses exact rational arithmetic and Python's standard library.

The domain calculation is part of the proof. The strongest certificate passes the complete $\Z[1/2]$ audit. The universal certificate and both standard-form certificates pass complete $\Z$ audits. In each case, the ring and rational closures agree on every queried tuple, establishing the bound in its stated group class.

Independently of our verifiers, companion certificates in the format of~\cite{Tomar} pass its public checker, run unchanged, at $1.8271149232\ldots$ (no $2$-torsion) and $1.8271176725\ldots$ (every abelian group), both below the previous record for the real constant.

\appendix
\section{Root and transported cross forms}\label{app:forms}
Let $f_0,\ldots,f_{42}$ be the root-page cross forms of~\cite{Tomar}, in the order of the supplied certificates, and let $f_i^{(j)}$ denote the same literal coefficient vector on page $P_j$. The complete root-form list is supplied in the certificate files and \texttt{anc/CERTIFICATE\_INDEX.md}. Table~\ref{tab:transport} gives every additional form on the copied pages. JSON indices are zero-based and refer to the displayed order within each certificate.
\begin{table}[H]
\centering\small
\begin{tabular}{@{}rl@{\quad}rl@{\qquad}rl@{\quad}rl@{}}\toprule
\multicolumn{4}{c}{Record: Theorem~\ref{thm:record}} & \multicolumn{4}{c}{All groups: Theorem~\ref{thm:all}}\\
\cmidrule(lr){1-4}\cmidrule(lr){5-8}
Index & Form & Index & Form & Index & Form & Index & Form\\\midrule
43 & $f_{1}^{(1)}$ & 64 & $f_{2}^{(5)}$ & 43 & $f_{29}^{(3)}$ & 64 & $f_{8}^{(3)}$\\
44 & $f_{2}^{(1)}$ & 65 & $f_{28}^{(3)}$ & 44 & $f_{11}^{(4)}$ & 65 & $f_{39}^{(6)}$\\
45 & $f_{13}^{(1)}$ & 66 & $f_{6}^{(3)}$ & 45 & $f_{13}^{(4)}$ & 66 & $f_{2}^{(5)}$\\
46 & $f_{21}^{(1)}$ & 67 & $f_{12}^{(4)}$ & 46 & $f_{2}^{(3)}$ & 67 & $f_{12}^{(3)}$\\
47 & $f_{28}^{(1)}$ & 68 & $f_{14}^{(4)}$ & 47 & $f_{40}^{(3)}$ & 68 & $f_{15}^{(3)}$\\
48 & $f_{29}^{(1)}$ & 69 & $f_{38}^{(5)}$ & 48 & $f_{28}^{(5)}$ & 69 & $f_{26}^{(5)}$\\
49 & $f_{40}^{(1)}$ & 70 & $f_{6}^{(4)}$ & 49 & $f_{36}^{(5)}$ & 70 & $f_{37}^{(5)}$\\
50 & $f_{29}^{(4)}$ & 71 & $f_{8}^{(4)}$ & 50 & $f_{29}^{(5)}$ & 71 & $f_{38}^{(5)}$\\
51 & $f_{28}^{(5)}$ & 72 & $f_{39}^{(1)}$ & 51 & $f_{40}^{(5)}$ & 72 & $f_{22}^{(6)}$\\
52 & $f_{11}^{(6)}$ & 73 & $f_{20}^{(1)}$ & 52 & $f_{2}^{(6)}$ & 73 & $f_{16}^{(5)}$\\
53 & $f_{13}^{(6)}$ & 74 & $f_{25}^{(3)}$ & 53 & $f_{13}^{(6)}$ & 74 & $f_{17}^{(5)}$\\
54 & $f_{2}^{(4)}$ & 75 & $f_{0}^{(4)}$ & 54 & $f_{21}^{(6)}$ & 75 & $f_{20}^{(6)}$\\
55 & $f_{40}^{(4)}$ & 76 & $f_{26}^{(1)}$ & 55 & $f_{40}^{(6)}$ & 76 & $f_{25}^{(2)}$\\
56 & $f_{36}^{(5)}$ & 77 & $f_{15}^{(4)}$ & 56 & $f_{1}^{(6)}$ & 77 & $f_{28}^{(2)}$\\
57 & $f_{29}^{(5)}$ &  &  & 57 & $f_{28}^{(6)}$ & 78 & $f_{6}^{(3)}$\\
58 & $f_{40}^{(5)}$ &  &  & 58 & $f_{29}^{(6)}$ & 79 & $f_{0}^{(3)}$\\
59 & $f_{31}^{(6)}$ &  &  & 59 & $f_{36}^{(6)}$ & 80 & $f_{15}^{(5)}$\\
60 & $f_{0}^{(1)}$ &  &  & 60 & $f_{1}^{(1)}$ & 81 & $f_{26}^{(2)}$\\
61 & $f_{22}^{(1)}$ &  &  & 61 & $f_{21}^{(1)}$ & 82 & $f_{0}^{(6)}$\\
62 & $f_{16}^{(5)}$ &  &  & 62 & $f_{28}^{(1)}$ & 83 & $f_{26}^{(6)}$\\
63 & $f_{17}^{(5)}$ &  &  & 63 & $f_{14}^{(3)}$ & 84 & $f_{31}^{(4)}$\\
\bottomrule\end{tabular}
\caption{Every copied-page cross form. Within each theorem, the indices run down the first pair of columns, then the second. Page superscripts follow Table~\ref{tab:glues}; form subscripts follow the root-form order defined above.}
\label{tab:transport}
\end{table}

\section{Certificates and data availability}\label{app:reproduction}
\begin{table}[H]
\centering\small
\begin{tabular}{@{}lrrrr@{}}\toprule
 & Record & All groups & $2762/1509$ & $152/83$\\\midrule
Vertices &42&42&108&30\\
Copy steps &9&9&20&7\\
Extra scalar forms &78&85&0&0\\
Relation rank &71&71&177&46\\
Quotient dimension &13&13&39&14\\
Submodularity rows &519&566&758&103\\
Monotonicity rows &0&1&0&0\\
Copy-identity rows &9&9&20&8\\
Lemma-application rows &12&11&0&0\\
Equal-law entries &167&179&568&59\\
Verified library lemmas &72&77&0&0\\
Domain-audit queries &1647&1779&2420&341\\\bottomrule
\end{tabular}
\caption{Sizes of the four exact elemental certificates.}
\label{tab:counts}
\end{table}

Four elemental certificates are supplied in \path{anc/certificates/}: \path{record.json}, \path{all-groups.json}, \path{standard.json}, and \path{standard-five-page.json}. The fields \texttt{steps}, \texttt{extra}, \texttt{eqlaw}, \texttt{rows} and \texttt{lambda} specify the conditional-copy configuration, extra forms, equal-law generators, weighted proof rows and target weights, respectively. A standard form is encoded by a vertex and its integer slope; \texttt{["amb",k]} denotes extra form $k$. Rational coefficients are stored as exact strings, and the verifier reconstructs the complete proof from these fields. The glue-word order in \eqref{eq:rootorder} differs from the JSON vertex numbering, which remains $(0,1,2,3,4,5)$.

Run the following from the source directory:
\begin{quote}\begin{minipage}{\linewidth}\ttfamily\small
python3 anc/verify\_group.py anc/certificates/record.json -{}-invert 2\\
python3 anc/verify\_group.py anc/certificates/all-groups.json\\
python3 anc/verify\_group.py anc/certificates/standard.json\\
python3 anc/verify\_group.py anc/certificates/standard-five-page.json
\end{minipage}\end{quote}
Each command first verifies the complete rational identity and lemma library using \texttt{stackcheck.py}, then performs the independent integer or localized-integer domain audit. Successful output is a JSON receipt with \texttt{PASS}, the exact constant and a passing record for every copy step. All four audits verify equality of the ring and rational closures on every queried tuple, with the counts in Table~\ref{tab:counts}. A failure returns nonzero status. The receipts bind the certificate bytes by SHA-256; \texttt{SHA256SUMS} binds the full distributed source and code. Run without Python's \texttt{-O} option.

For the second rational replay of either record, run
\begin{quote}\begin{minipage}{\linewidth}\ttfamily\small
python3 anc/check630\_extended.py anc/certificates/record.json\\
python3 anc/check630\_extended.py anc/certificates/all-groups.json
\end{minipage}\end{quote}
This wrapper retains the Shannon-row and auxiliary-lemma verification code of~\cite{Tomar}, uses sparse exact rational closure, and supports the shared and transported cross-form equal laws of Lemma~\ref{lem:enriched}. The two simpler certificates replay with \texttt{check630\_fast.py}. The original \texttt{check630.py} is included unmodified under its Apache-2.0 licence.

\begin{samepage}
Two companion certificates have completed replays with the previous record's public checker~\cite{Tomar}, run unchanged. The file \texttt{record-original-format.json} proves
\[
\frac{2858285574867171691166}{1564370986493012722265}
=1.8271149232\ldots
\]
in groups without $2$-torsion; \texttt{all-groups-original-format.json} proves
\[
\frac{226483966351349932879}{123956967720824699832}
=1.8271176725\ldots
\]
in every abelian group. Each receipt binds its own companion input and the unchanged checker by SHA-256. The supplement's \texttt{CERTIFICATE\_INDEX.md} tabulates the exact exponents, weight triples, proof sizes, hashes and passing checker receipts; \texttt{SHA256SUMS} binds the complete release.
\end{samepage}

\subsection{Controls and provenance}
The driver \texttt{negative\_controls.py} deliberately alters target weights, Shannon multipliers, copy coefficients, equal-law entries, glue data and auxiliary lemma data. All 34 deliberately corrupted inputs were rejected. Four additional corruptions of cross-form support and images were also rejected by the extended replay. Five further controls, each corrupting a coefficient in one of the newly added Dougherty--Freiling--Zeger proofs, were rejected as well. The supplement records the complete results and the exact inputs tested.

The source package contains the manuscript, complete certificate data, checking programs, licences and observed verification receipts. The file \texttt{SHA256SUMS} binds all distributed proof and source files.

\section*{Research methods}
CASPAC, the computational mathematics engine of CASP, discovered the certificates through AI-directed search under the author's direction. Independent exact checkers verified the resulting proofs.

\bigskip\noindent{\small\textsc{Atrin Farnamfar}\quad\texttt{atrinfarnamfar@gmail.com}}
\end{document}